\documentclass[11pt,a4paper,reqno]{amsart}

\usepackage{amsmath,amssymb,amsthm,mathtools}
\usepackage[margin=1in]{geometry}
\usepackage{iftex}
\ifPDFTeX
 \usepackage[T1]{fontenc}
 \usepackage{lmodern}
\fi
\usepackage[expansion=false]{microtype}
\usepackage{tikz}
\usepackage[hidelinks]{hyperref}
\hypersetup{
 pdftitle={A seven-element interval of power-semiring varieties},
 pdfauthor={Aifa Wang, Lili Wang, Qingrui Yin},
 pdfsubject={Intervals of additively idempotent semiring varieties}
}

\newtheorem{theorem}{Theorem}[section]
\newtheorem{proposition}[theorem]{Proposition}
\newtheorem{lemma}[theorem]{Lemma}
\newtheorem{corollary}[theorem]{Corollary}
\theoremstyle{remark}
\newtheorem{remark}[theorem]{Remark}
\numberwithin{equation}{section}

\newcommand{\AI}{\mathsf{AI}}
\newcommand{\V}{\mathsf V}
\newcommand{\Pplus}{\mathcal P^+}
\newcommand{\Mod}{\operatorname{Mod}}
\newcommand{\Dvar}{\mathsf D}
\newcommand{\B}{\mathcal B}
\newcommand{\Rvar}{\mathsf R}

\title[An interval of power-semiring varieties]{A seven-element interval of power-semiring varieties}
\author[A. Wang]{Aifa Wang}
\address[Aifa Wang]{School of Mathematical Sciences,\newline Chongqing University of Technology, Chongqing 400054, P.R. China}
\email{wangaf@cqut.edu.cn}
\thanks{Corresponding author: Aifa Wang.}
\author[L. Wang]{Lili Wang}
\address[Lili Wang]{School of Mathematical Sciences,\newline Chongqing University of Technology, Chongqing 400054, P.R. China}
\email{wllaf@cqut.edu.cn}
\author[Q. Yin]{Qingrui Yin}
\address[Qingrui Yin]{School of Mathematical Sciences,\newline Chongqing University of Technology, Chongqing 400054, P.R. China}
\email{13996559059@163.com}
\date{}

\begin{document}

\begin{abstract}
We determine the interval between the varieties generated by the nonempty
and full power semirings of a semigroup satisfying \(xy\approx x^2\).
It has two elements for a one-element semigroup, three for a nontrivial
left-zero or null semigroup, and seven otherwise, without any finiteness
assumption.  In the seven-element case the interval is
\([\Rvar,\Rvar^0]\), where
\(\Rvar=\Mod_{\rm ai}\{xy\approx xz\}\) and \(\Rvar^0\) is its
zero-adjunction variety.  Every word absorption valid in \(\Rvar\) but
not in \(\Rvar^0\) is equivalent over \(\Rvar^0\) to one of six
relative identities.  We give exact absorption criteria and finite bases
for all seven varieties.  Three have basis rank two and four have basis
rank three; their two-generated free algebras have orders \(15\), \(24\),
or \(35\).  We also count the finite-rank free algebras of both endpoints.
\end{abstract}


\maketitle

\section{Introduction}

We work with additively idempotent semirings in the signature
\((+,\cdot)\), without named constants.  For a semigroup \(S\), let
\(\Pplus(S)\) and \(\mathcal P(S)\) denote its nonempty and full power
semirings, with union and setwise multiplication.  Thus
\(\mathcal P(S)\) is obtained from \(\Pplus(S)\) by adjoining the empty
set as an additively least multiplicative zero.  Put
\[
 I(S)=[\V(\Pplus(S)),\V(\mathcal P(S))].
\]
Yue, Ren, and Gao~\cite[Problem~1.3]{YueRenGao2026} ask for the
cardinality of \(I(S)\).  They prove that it is \(2^{\aleph_0}\) when
\(S\) is the multiplicative reduct of
\(S_7\)~\cite[Theorem~5.2]{YueRenGao2026}.

We determine \(I(S)\) for every semigroup satisfying
\begin{equation}\label{eq:semigroup-identity}
xy\approx x^2.
\end{equation}
Such a semigroup has multiplication \(xy=f(x)\), where \(f^2=f\).
The interval has two elements if \(|S|=1\), three if \(S\) is a
nontrivial left-zero or null semigroup, and seven if
\(f\ne\operatorname{id}_S\) and \(|f(S)|\geq2\).
In the last case both endpoint varieties are generated by the power
semirings of the three-element subsemigroup
\begin{equation}\label{eq:three-element-table}
\begin{array}{c|ccc}
\cdot&a&b&c\\ \hline
a&a&a&a\\
b&b&b&b\\
c&a&a&a
\end{array}.
\end{equation}

The lower endpoint in this case is
\[
 \Rvar=\Mod_{\rm ai}\{xy\approx xz\}.
\]
Yue, Ren, Zeng, and Shao established its finite
generation~\cite[Proposition~4.5 and Remark~4.6]{YueRenZengShao2025};
Yue and Ren determined its ten-element subvariety lattice and proved
that it is a Cross
variety~\cite[Lemma~2.1, Theorem~2.7 and Corollary~2.8]{YueRen2025}.
Here we determine the interval \([\Rvar,\Rvar^0]\), where
\(\Rvar^0\) is generated by the zero extensions of members of \(\Rvar\).

The proof classifies the relative identities of these endpoints.
Every word absorption valid in \(\Rvar\) but not in \(\Rvar^0\)
reduces, over \(\Rvar^0\), to one of six identities.
Exact absorption criteria distinguish the resulting varieties and
establish their inclusion relations.  The criteria follow from word
normal forms and invariants under substitutions and term contexts;
explicit substitutions give the six-type reduction.
All seven varieties are finitely based.

We also determine their basis ranks.  The seven varieties have only
three distinct two-variable equational theories, with two-generated
free algebras of orders \(15\), \(24\), and \(35\).
Three varieties are generated by these free algebras; the other four
require three free generators.  Counting formulas for the endpoint
free algebras follow from their normal forms.

The upper-endpoint absorption criterion is a specialization of the
zero-adjunction criterion of Wu, Ren, and
Zhao~\cite[Proposition~1.5]{WuRenZhao2024}, restated in
Gao, Ren, and Zhao~\cite[Proposition~2.1]{GaoRenZhao2026}.
The latter paper also shows that zero adjunction preserves equality
of generated varieties.  Together with the known four-element
generator of \(\Rvar\), this gives a five-element generator of
\(\Rvar^0\) (Proposition~\ref{prop:zero-endpoint}).
For finite groups of order at least three, Shao and
Gao~\cite[Corollary~1.2]{ShaoGao2026} prove that every variety in the
corresponding power-semiring interval has infinite axiomatic rank and
is nonfinitely based; their \(\mathcal P(G)\) denotes the nonempty
power semiring.
Related classifications for semilattice-ordered semigroups and normal
bands appear in Ku\v ril and Pol\'ak~\cite{KurilPolak2005} and
Ku\v ril~\cite{Kuril2017}.  Ren and Zhao~\cite{RenZhao2016} determine
a nine-element variety lattice under \(x^3\approx x\) and
\(xy\approx yx\).  For general variety theory we refer to
Burris and Sankappanavar~\cite{BurrisSankappanavar}.
\section{Preliminaries and the main theorem}

Let \(\AI\) denote the variety of additively idempotent semirings in the signature \((+,\cdot)\).  For terms \(u,v\), write
\[
 u\leq v \quad\text{for the identity}\quad u+v\approx v.
\]
Every term is, by distributivity and additive idempotence, a finite nonempty
sum of nonempty words.  We identify such a polynomial \(p\) with its finite
set of word summands and write \(w\in p\) for membership in that set.
For a word \(w\), let \(|w|\), \(h(w)\), and \(c(w)\) denote its length,
first variable, and content, respectively.  A word of length one is called
a \emph{one-letter word}; a word of length at least two is called a
\emph{long word}.  For \(h\in C\), where \(C\) is a finite nonempty variable
set, fix a long word \(w_{h,C}\) with first variable \(h\) and content \(C\).
For \(C=\{h\}\), take \(w_{h,C}=h^2\); otherwise take \(h\) followed by
the remaining elements of \(C\) in some fixed order.
For a polynomial \(p\), set
\begin{equation}\label{eq:LE}
\begin{aligned}
 L(p)&=\{x:x\text{ occurs as a one-letter summand of }p\},\\
 E(p)&=\{h(w):w\text{ is a long summand of }p\}.
\end{aligned}
\end{equation}

\begin{theorem}[classification]\label{thm:main}
Let \(S\neq\varnothing\) be a semigroup satisfying \(xy\approx x^2\).  Then there is an idempotent map \(f:S\to S\) such that \(xy=f(x)\), and
\begin{equation}\label{eq:main-cardinality}
|I(S)|=
 \begin{cases}
 2,& |S|=1,\\
 3,& |S|\geq2\text{ and }(f=\operatorname{id}_S\text{ or }|f(S)|=1),\\
 7,& f\neq\operatorname{id}_S\text{ and }|f(S)|\geq2.
 \end{cases}
\end{equation}
In the last case the interval is, up to isomorphism, the seven-element lattice described in Theorem~\ref{thm:lattice} below.
\end{theorem}

\begin{proof}
If \(xy=f(x)\) with \(f^2=f\), then
\[
 (xy)z=f(f(x))=f(x)=x(yz),
\]
so the operation is associative and \(x^2=f(x)\).  Conversely, if
\(xy\approx x^2\), define \(f(x)=x^2\).  Taking \(y=x^2\) gives
\(x^3=x^2\), whence \(x^4=x^2\) and \(f^2=f\).
This proves the structural assertion.  The nondegenerate case of the
cardinality assertion is proved in Theorem~\ref{thm:lattice}; all remaining
cases are proved in Section~\ref{sec:boundary}.
\end{proof}

For the remainder of the paper put
\begin{equation}\label{eq:endpoints}
A=\V(\Pplus(S)),\qquad U=\V(\mathcal P(S)).
\end{equation}
The symbol \(\V\) denotes the variety generated by an algebra, and
\(\Mod_{\rm ai}\Sigma\) denotes the class of algebras in \(\AI\)
satisfying \(\Sigma\).  All inclusions of varieties are in this fixed
signature.  The cardinality of an interval includes its endpoints.
For an ai-semiring \(T\), let \(T^0\) be obtained by adjoining a new
additively least multiplicative zero, with the original operations on
\(T\) unchanged.  For a variety \(W\), put
\[
 W^0=\V\{T^0:T\in W\}.
\]
The zero is an element of these algebras, not a named constant in the
signature.  In particular,
\(\mathcal P(S)\cong(\Pplus(S))^0\).

\section{Endpoint theories and a basis for the upper endpoint}\label{sec:endpoints}

Assume throughout this section and Sections~\ref{sec:relative}--\ref{sec:seven} that
\begin{equation}\label{eq:nondegeneracy}
f\neq\operatorname{id}_S,\qquad |f(S)|\geq2.
\end{equation}
Since an idempotent surjective map is the identity, choose \(c\notin f(S)\), put \(a=f(c)\), and choose \(b\in f(S)\setminus\{a\}\).  The set \(S_*:=\{a,b,c\}\) is a subsemigroup, with multiplication table \eqref{eq:three-element-table}.  In particular, the endpoint theories of every such \(S\) are already visible on this three-element subsemigroup.

\begin{lemma}[endpoint absorption criterion]\label{lem:endpoint}
For every polynomial \(p\),
\begin{equation}\label{eq:one-letter-criterion}
A\models x\leq p \quad\Longleftrightarrow\quad U\models x\leq p
 \quad\Longleftrightarrow\quad x\in L(p).
\end{equation}
If \(u=w_{h,C}\) is a long word with first variable \(h\) and content \(C\), then
\begin{equation}\label{eq:lower-criterion}
A\models u\leq p\quad\Longleftrightarrow\quad h\in E(p),
\end{equation}
whereas
\begin{equation}\label{eq:upper-criterion}
U\models u\leq p\quad\Longleftrightarrow\quad
 \exists v\in p\;\bigl(|v|\geq2,\ h(v)=h,\ c(v)\subseteq C\bigr).
\end{equation}
Consequently \(\V(\Pplus(S))=\V(\Pplus(S_*))\) and \(\V(\mathcal P(S))=\V(\mathcal P(S_*))\).
\end{lemma}

\begin{proof}
Write \(f(X)=\{f(s):s\in X\}\) for the direct image of a subset.
For nonempty variable sets \(X_y\subseteq S\), every long word \(w\)
evaluates in \(\Pplus(S)\) to
\[
 w(X_1,\ldots,X_n)=f(X_{h(w)}).
\]
In the full power semiring the same formula holds when
\(X_y\neq\varnothing\) for every \(y\in c(w)\).
Its value is empty if \(X_y=\varnothing\) for some \(y\in c(w)\).
This proves the sufficiency of \eqref{eq:one-letter-criterion}--\eqref{eq:upper-criterion}.

For necessity it is enough to use the three-element table \eqref{eq:three-element-table}.
If \(x\notin L(p)\), assign \(X_x=\{c\}\) and every other variable
the value \(\{b\}\).  Every right-hand long summand has value
\(\{a\}\) or \(\{b\}\), and every right-hand one-letter summand has value
\(\{b\}\).  Thus \(p\) does not contain \(c\), whereas \(x\) does.
If \(h\notin E(p)\), assign \(X_h=\{c\}\) and all other variables
\(\{b\}\).  The long left-hand word has value \(\{a\}\), while the
right-hand side is contained in \(\{b,c\}\).
Finally, if the condition in \eqref{eq:upper-criterion} fails, assign variables outside
\(C\) the empty set, assign \(h\) the set \(\{c\}\), and assign the
remaining variables of \(C\) the set \(\{b\}\).
The left-hand word has value \(\{a\}\); surviving right-hand long words
cannot start with \(h\), so no right-hand summand contains \(a\).
Every identity is equivalent to word absorptions in both directions.
The criteria therefore describe the entire endpoint theories, both for
\(S\) and for \(S_*\), proving the final assertion.
\end{proof}

\begin{proposition}[finite basis for \(U\)]\label{prop:upperbasis}
Relative to the axioms of \(\AI\), the upper endpoint \(U\) has the basis
\begin{equation}\label{eq:upper-basis}
\B=\{xy\approx x^2y,\quad xy\approx xy^2,\quad
 xyz\approx xzy,\quad xyz\leq xy\}.
\end{equation}
\end{proposition}

\begin{proof}
The four identities hold in \(\mathcal P(S)\) by the endpoint computation.
The identity \(xyz\approx xzy\), applied to prefixes and in right contexts,
permutes any adjacent pair of letters after the first position.
The identity \(xy\approx xy^2\) removes repetitions in those positions.
If a word contains a variable other than its first variable, tail
permutation followed by \(x^2y\approx xy\) removes extra occurrences of
the first variable.  A word whose content is \(\{h\}\) reduces to
\(h^2\), since \(x^2\approx x^3\) follows from the first identity.
Thus every long word reduces to \(w_{h,C}\).

Write a long word \(v\) as \(v=v_1v_2\) with both factors nonempty.
Substituting \(x=v_1,y=v_2,z=t\) in the last identity gives
\(vt\leq v\).  If long words \(u,v\) have the same first variable and
\(c(v)\subseteq c(u)\), append the missing variables to \(v\) and
normalize to obtain \(u\leq v\).
Lemma~\ref{lem:endpoint} now derives every word absorption valid in \(U\).
An identity between sums is equivalent to absorbing every summand of each
side into the other side, proving completeness.
\end{proof}

We shall repeatedly use two consequences:
every long word \(u\) satisfies \(u^2\approx u\) in \(U\), and
\(uv\approx u\) whenever \(c(v)\subseteq c(u)\).
These follow by the same normal forms.

\begin{remark}[the lower endpoint]\label{rem:lower}
The lower endpoint in the nondegenerate case is the already studied variety
\begin{equation}\label{eq:known-lower}
A=\operatorname{Mod}_{\rm ai}\{xy\approx xz\}.
\end{equation}
Indeed, \(xy\approx xz\) implies \(xy\approx x^2\) by putting \(z=x\), and then every long word reduces to the square of its first variable.  The criterion \eqref{eq:one-letter-criterion}--\eqref{eq:lower-criterion} shows that no additional identity identifies distinct \(L\)- or \(E\)-data.  Thus \eqref{eq:known-lower} follows.  This endpoint is the variety \(\Rvar\) (denoted by \(\mathbf R\)) studied by Yue and Ren~\cite{YueRen2025}; its finite generation is due to Yue, Ren, Zeng, and Shao~\cite{YueRenZengShao2025}, and its ten-element subvariety lattice and Cross property are given in~\cite[Theorem~2.7 and Corollary~2.8]{YueRen2025}.
\end{remark}

\begin{proposition}[zero-adjunction description and a small generator]
\label{prop:zero-endpoint}
Under \eqref{eq:nondegeneracy},
\[
 A=\Rvar,\qquad U=\Rvar^0,\qquad
 I(S)=[\Rvar,\Rvar^0].
\]
If \(K=S_{(4,475)}\) is the four-element ai-semiring in
\cite[Table~4 and Lemma~2.1]{YueRen2025}, then
\[
 A=\V(K),\qquad U=\V(K^0).
\]
In particular, the upper endpoint is generated by a five-element algebra.
\end{proposition}

\begin{proof}
Remark~\ref{rem:lower} gives \(A=\Rvar\), and
\cite[Lemma~2.1]{YueRen2025} gives \(\Rvar=\V(K)\).
By the general zero-adjunction result
\cite[Corollaries~2.2 and~2.3]{GaoRenZhao2026},
\(\V(T)^0=\V(T^0)\) for every ai-semiring \(T\).
Applying it to \(T=\Pplus(S)\) and to \(T=K\) yields
\[
 U=\V((\Pplus(S))^0)=A^0=\Rvar^0=\V(K^0).
\]
Since zero adjunction adds one new element, \(|K^0|=5\).
\end{proof}

\begin{remark}
The five-element generator in Proposition~\ref{prop:zero-endpoint}
is not asserted to be free.  The orders of free algebras and the least
number of free generators needed to generate each interval variety
are determined in Sections~\ref{sec:seven} and~\ref{sec:ranks}.
\end{remark}

\section{Six relative absorption rules}\label{sec:relative}

Define the following six identities, all written as absorptions:
\begin{equation}\label{eq:relative-rules}
\begin{aligned}
 \alpha &: xy\leq x+y^2+xz, &\qquad
 \beta &: xy\leq y^2+xz,\\
 \eta &: x^2\leq x+xy, &
 \delta &: xy\leq y+xz+yz,\\
 \gamma &: xy\leq y+xz, &
 \sigma &: x^2\leq xy.
\end{aligned}
\end{equation}
Let \(U_\theta=U\cap\Mod\{\theta\}\) for \(\theta\) among these six rules.  Every rule in \eqref{eq:relative-rules} holds in \(A\).

Fix a finite nonempty set \(C\) of variables and \(h\in C\).  For a polynomial \(p\), define
\begin{equation}\label{eq:predicates}
\begin{aligned}
 H_{h,C}(p)&\Longleftrightarrow\exists v\in p\ (h(v)=h,\ c(v)\subseteq C),\\
 N_C(p)&\Longleftrightarrow\exists v\in p\ (|v|\geq2,\ c(v)\subseteq C),\\
 T_C(p)&\Longleftrightarrow C\cap L(p)\cap E(p)\neq\varnothing,\\
 J_C(p)&\Longleftrightarrow\exists v\in p\ (c(v)\subseteq C),\\
 Q_C(p)&\Longleftrightarrow N_C(p)\lor T_C(p).
\end{aligned}
\end{equation}
Thus \(H\) and \(J\) may be witnessed by a one-letter word, while \(N\) must be witnessed by a long word.

\begin{lemma}[contextual invariants]\label{lem:invariants}
Every identity valid in \(U\) preserves the seven data \(L,E,H,N,T,J,Q\).  After arbitrary polynomial substitutions and insertion into an arbitrary one-hole term context, the two sides of one of the rules in \eqref{eq:relative-rules} preserve the data in the following table:
\begin{equation}\label{eq:invariants}
\begin{array}{c|c}
\text{rule}&\text{data preserved for every }(h,C)\\ \hline
\alpha&L,E,H,N\\
\beta&L,E,N\\
\eta&L,E,H\\
\delta&L,E,Q\\
\gamma&L,E,J\\
\sigma&L,E.
\end{array}
\end{equation}
\end{lemma}

\begin{proof}
If \(p\approx q\) holds in \(U\), \eqref{eq:one-letter-criterion}--\eqref{eq:upper-criterion} say that the
one-letter summands agree and each long summand of either side is
covered by a long summand on the other side with the same first variable
and smaller content.  Thus \(L,E,H,N,J\) agree, and so do \(T,Q\).
This applies in particular to any substituted and contextual instance of
an identity of \(U\).

Consider a rule \(l\leq r\), regarded as \(l+r\approx r\), and a
substitution \(\varphi\).  Expanding a one-hole context by distributivity
expresses its polynomial as a common part together with finitely many
pieces \(\lambda\varphi(l+r)\rho\) or \(\lambda\varphi(r)\rho\).
Here \(\lambda,\rho\) are possibly empty words; this notation does not
introduce constants into the signature.  Since the second expansion is
contained in the first, it suffices to transfer every witness supplied
by a deleted summand to a summand common to both sides.

For \(\alpha,\beta,\gamma,\delta\), a deleted summand has the form
\(\lambda ab\rho\), with \(a\in\varphi(x)\), \(b\in\varphi(y)\).
For any \(d\in\varphi(z)\), the common expansion contains
\(\lambda ad\rho\).  These are long words with the same first variable.
For \(\eta,\sigma\), a deleted summand has the same form, but now
\(a,b\in\varphi(x)\) are chosen independently; choose
\(d\in\varphi(y)\) to get the same protection.
All substituted polynomials are nonempty, so \(d\) exists.
Thus all six rules preserve \(E\), and they preserve \(L\) because every
deleted summand is long.

For \(\alpha\), the common summand \(\lambda a\rho\) has the same
first variable and smaller content than \(\lambda ab\rho\), so it
protects \(H\).  The common long summand \(\lambda b^2\rho\) protects
\(N\), since its content is contained in that of the deleted summand.
This also proves \(N\)-preservation for \(\beta\).
For \(\eta\), \(\lambda a\rho\) protects \(H\), even when the two
selected words \(a,b\in\varphi(x)\) differ.
For \(\gamma\), \(\lambda b\rho\) protects \(J\).

For \(\delta\), suppose the deleted summand witnesses \(N_C\).
The common summand \(\lambda b\rho\) also has content in \(C\).
If it is long, it witnesses \(N_C\).  Otherwise
\(\lambda=\rho=\varepsilon\) and \(b=k\) is a variable in \(C\).
The common right-hand side contains both \(k\), from \(\varphi(y)\),
and \(kd\), from \(\varphi(yz)\).
Hence \(k\in C\cap L\cap E\), which witnesses \(T_C\).
Already existing \(T_C\)-witnesses are preserved because \(L,E\) are.
Therefore \(Q_C=N_C\lor T_C\) is preserved.

Applying these transfers to every deleted summand proves equality of the
specified data on the two sides of each elementary replacement.
It applies to replacements in either direction.
Induction on an equational derivation proves the assertion.
\end{proof}

\begin{theorem}[exact relative criteria]\label{thm:relativecriteria}
For every \(\theta\) in \eqref{eq:relative-rules}, absorption of a one-letter word in
\(U_\theta\) is characterized by \eqref{eq:one-letter-criterion}.
For a long word \(u=w_{h,C}\), the necessary and sufficient conditions for
\(U_\theta\models u\leq p\) are
\begin{equation}\label{eq:relative-criteria}
\begin{array}{c|c}
\text{variety}&\text{condition}\\ \hline
U_\alpha&h\in E(p),\ H_{h,C}(p),\ N_C(p)\\
U_\beta&h\in E(p),\ N_C(p)\\
U_\eta&h\in E(p),\ H_{h,C}(p)\\
U_\delta&h\in E(p),\ N_C(p)\lor T_C(p)\\
U_\gamma&h\in E(p),\ J_C(p)\\
U_\sigma&h\in E(p).
\end{array}
\end{equation}
In particular, \(U_\sigma=A\).
\end{theorem}

\begin{proof}
The necessity of \(x\in L(p)\) is \eqref{eq:one-letter-criterion}, and the necessity of \(h\in E(p)\) is \eqref{eq:lower-criterion}.  The remaining necessities follow from Lemma~\ref{lem:invariants}: the left-hand word itself witnesses the relevant data, and every derivation from the defining rule must preserve that data.  Equational completeness turns this non-derivability statement into failure in the corresponding variety.

For sufficiency, choose a long summand \(g\in p\) with first variable \(h\), and choose an auxiliary word \(t\) whose content contains all variables needed to ensure \(ht\leq g\).  For \(\sigma\), substitute \(x=u,y=t\) and use \(u^2\approx u\) to obtain \(u\leq ut\leq p\).  If \(J_C(p)\) holds, choose a witness \(v\) and use \(x=u,y=v,z=t\) in \(\gamma\); since \(c(v)\subseteq C\), \(uv\approx u\), and hence \(u\leq v+ut\leq p\).  The same substitution with an internal long witness and \(v^2\approx v\) proves the \(\beta\)-case.

For \(\eta\), if \(H\) is witnessed by a long word, \eqref{eq:upper-criterion} already applies.  Otherwise \(h\in L(p)\); with \(x=h,y=t\), the identity \(x^2\leq x+xy\), together with \(u\leq h^2\) and \(ht\leq g\), gives the result.  For \(\alpha\), the only additional case is again \(h\in L(p)\).  Choose an internal long witness \(v\) for \(N\), use \(x=h,y=v,z=t\), and obtain
\[
 u\leq hv\leq h+v^2+ht\approx h+v+ht\leq p.
\]
Finally suppose that \(\delta\) applies.  If \(N\) holds, use an internal long witness \(v\), substitute \(x=u,y=v,z=t\), and use \(vt\leq v\) to get
\[
 u\approx uv\leq v+ut+vt\leq p.
\]
If \(T\) holds, choose \(k\in C\cap L(p)\cap E(p)\), and choose a long right-hand word \(g_k\) with first variable \(k\).  Taking \(t\) to contain the contents of both \(g\) and \(g_k\), the substitution \(x=u,y=k,z=t\) gives
\[
 u\approx uk\leq k+ut+kt\leq p.
\]
This proves all rows of \eqref{eq:relative-criteria}.  The last row is exactly the endpoint criterion \eqref{eq:lower-criterion}, so \(U_\sigma=A\).
\end{proof}

\section{The six-type reduction}\label{sec:reduction}

\begin{theorem}[six types]\label{thm:sixtypes}
Let \(u\) be a word and \(p\) a polynomial.  Suppose that \(A\models u\leq p\) but \(U\not\models u\leq p\).  Then the absorption is equivalent, over \(U\), to exactly one of the six rules in \eqref{eq:relative-rules}.  More precisely, with \(u=w_{h,C}\), the six mutually exclusive cases are
\begin{equation}\label{eq:six-types}
\begin{array}{c|c}
\text{condition}&\text{equivalent rule}\\ \hline
H\land N&\alpha\\
H\land\neg N&\eta\\
\neg H\land N&\beta\\
\neg H\land\neg N\land T&\delta\\
\neg H\land\neg N\land\neg T\land J&\gamma\\
\neg J&\sigma.
\end{array}
\end{equation}
Here all predicates refer to \((h,C,p)\).
\end{theorem}

\begin{proof}
By Lemma~\ref{lem:endpoint}, \(u\) is long and \(h\in E(p)\).
Call a word \emph{internal} if its content is contained in \(C\), and
\emph{external} otherwise.  The failure of \eqref{eq:upper-criterion} says that no internal
long right-hand word starts with \(h\).  Hence, in the strict situation,
\(H\) is equivalent to \(h\in L(p)\).  In particular \(H\) implies
\(T\), and \(N\) or \(T\) implies \(J\), proving exhaustiveness and
mutual exclusivity of the cases in \eqref{eq:six-types}.

In each case below, a substitution \(\varphi\) gives
\(\varphi(u)\leq\varphi(p)\leq R\); the second absorption is derived
using only the identities of \(U\).

\smallskip\noindent
\emph{The cases with \(N\).}
An internal long right-hand word starts with a variable other than \(h\),
so \(|C|\geq2\).  Use
\[
 \varphi(k)=
 \begin{cases}
 x,&k=h,\\
 y^2,&k\in C\setminus\{h\},\\
 xz,&k\notin C.
 \end{cases}
\]
The left side normalizes to \(xy\).
Internal long words and internal one-letter words other than \(h\) are
absorbed by \(y^2\).
An external word whose image starts with \(x\) contains \(z\) and is
absorbed by \(xz\); if its image starts with \(y\), it is absorbed by
\(y^2\).  The possible one-letter summand \(h\) contributes \(x\).
Thus \(R=x+y^2+xz\) if \(H\) holds and \(R=y^2+xz\) otherwise.
These are \(\alpha\) and \(\beta\), respectively.

\smallskip\noindent
\emph{The case \(H\land\neg N\).}
Map internal variables to \(x\) and external variables to \(xy\).
The left side becomes \(x^2\).  Every internal right-hand word is
one-letter and becomes \(x\), while every external word becomes a long
word with first variable \(x\) and content \(\{x,y\}\).
Thus \(R=x+xy\), giving \(\eta\).

\smallskip\noindent
\emph{The case \(\neg H\land\neg N\land T\).}
Put \(M=C\cap L(p)\).  Here \(M\neq\varnothing\) and \(h\notin M\).
Use
\[
 \varphi(k)=
 \begin{cases}
 y,&k\in M,\\
 x,&k\in C\setminus M,\\
 xz,&k\notin C.
 \end{cases}
\]
The left side becomes \(xy\), since its content includes \(h\) and
some element of \(M\).  Internal right-hand words are one-letter
summands from \(M\), so they become \(y\).
External words whose first variable is in \(M\) are absorbed by \(yz\);
all other external words are absorbed by \(xz\).
Thus \(R=y+xz+yz\), giving \(\delta\).

\smallskip\noindent
\emph{The case \(\neg H\land\neg N\land\neg T\land J\).}
Since \(J\) holds and \(N\) fails, \(M=C\cap L(p)\) is again nonempty.
Use the preceding substitution.
Now \(M\cap E(p)=\varnothing\), since \(T\) fails.
No external long word starts in \(M\), and an external one-letter word
cannot lie in \(M\).
Consequently \(R=y+xz\), giving \(\gamma\).

\smallskip\noindent
\emph{The case \(\neg J\).}
Map internal variables to \(x\) and external variables to \(xy\).
Every right-hand word contains an external variable and normalizes to
\(xy\), whereas the left side normalizes to \(x^2\).
This gives \(\sigma\).

\smallskip
In each row, the indicated rule implies the original absorption by
Theorem~\ref{thm:relativecriteria}.  The rules are pairwise
inequivalent, as will follow from the separation argument in
Theorem~\ref{thm:lattice}.
\end{proof}

\section{The seven-element interval}\label{sec:seven}

\begin{figure}[htbp]
\centering
\begin{tikzpicture}[x=1.6cm,y=0.73cm,
 every node/.style={font=\small}]
 \coordinate (s) at (0,0);
 \coordinate (g) at (0,1);
 \coordinate (d) at (0,2);
 \coordinate (b) at (-0.75,3);
 \coordinate (e) at (0.75,3);
 \coordinate (a) at (0,4);
 \coordinate (u) at (0,5);
 \draw (s)--(g)--(d)--(b)--(a)--(u);
 \draw (d)--(e)--(a);
 \foreach \v in {s,g,d,b,e,a,u}
  \filldraw[fill=white] (\v) circle (2pt);
 \node[right=5pt] at (s) {\(A=U_\sigma\)};
 \node[right=5pt] at (g) {\(U_\gamma\)};
 \node[right=5pt] at (d) {\(U_\delta\)};
 \node[left=5pt] at (b) {\(U_\beta\)};
 \node[right=5pt] at (e) {\(U_\eta\)};
 \node[right=5pt] at (a) {\(U_\alpha\)};
 \node[right=5pt] at (u) {\(U\)};
\end{tikzpicture}
\caption{The seven-element interval, ordered by inclusion.}
\label{fig:interval}
\end{figure}
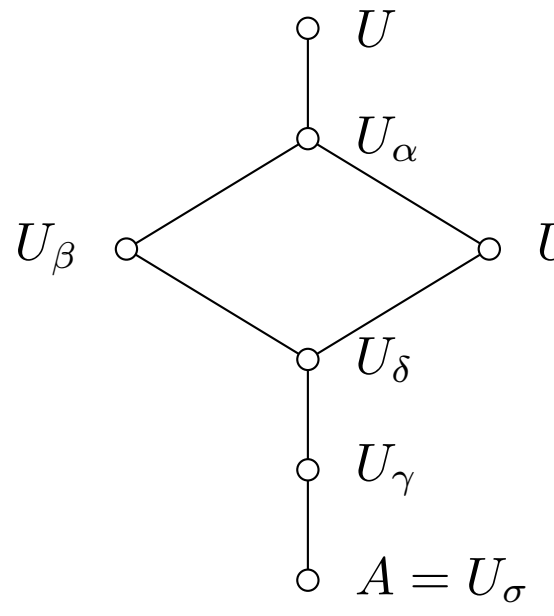

\begin{theorem}[the lattice]\label{thm:lattice}
Under the nondegeneracy assumption \eqref{eq:nondegeneracy},
\begin{equation}\label{eq:seven-varieties}
I(S)=\{U,U_\alpha,U_\beta,U_\eta,U_\delta,U_\gamma,U_\sigma=A\},
\end{equation}
and the covering relations are
\begin{equation}\label{eq:covers}
\begin{gathered}
 U\succ U_\alpha,\qquad U_\alpha\succ U_\beta,U_\eta,\\
 U_\beta,U_\eta\succ U_\delta,\qquad U_\delta\succ U_\gamma\succ A.
\end{gathered}
\end{equation}
The varieties \(U_\beta\) and \(U_\eta\) are incomparable.
The Hasse diagram is shown in Figure~\ref{fig:interval}.
\end{theorem}

\begin{proof}
The implications between the six rules are
\begin{equation}\label{eq:implications}
\sigma\Longrightarrow\gamma\Longrightarrow\delta,\qquad
 \delta\Longrightarrow\beta,\eta,\qquad
 \beta,\eta\Longrightarrow\alpha.
\end{equation}
For example, in \(\delta\) put \(y=x,z=y\) to obtain \(\eta\); put \(y=y^2\) and use the upper identities to obtain \(\beta\).  Conversely, \(\beta\) and \(\eta\) give
\[
 xy\leq y^2+xz\leq y+yz+xz,
\]
so together they are equivalent to \(\delta\).  Also \(\eta\) implies \(xy\leq x^2\leq x+xz\), hence \(\alpha\).

Let \(W\) be any variety in the interval.  An identity between sums is equivalent to its word absorptions.  Since \(A\subseteq W\), every absorption valid in \(W\) is valid in \(A\).  Remove the absorptions already valid in \(U\).  Theorem~\ref{thm:sixtypes} replaces each remaining absorption by one of the six rules, and \eqref{eq:implications} reduces every set of rules to no rule or to exactly one of the six displayed nodes.  Thus there are at most seven varieties.

All seven candidates are distinct.  The rule \(\alpha\) fails in \(U\)
under \(x=\{c\},y=\{b\},z=\varnothing\) in \eqref{eq:three-element-table}:
the left side contains \(a\) and the right side does not.
For the other separations, the predicates associated with the
left-hand word and right-hand polynomial of each rule are
\[
\begin{array}{c|cccc}
 &H&N&T&J\\ \hline
\alpha&1&1&1&1\\
\beta&0&1&0&1\\
\eta&1&0&1&1\\
\delta&0&0&1&1\\
\gamma&0&0&0&1\\
\sigma&0&0&0&0
\end{array}.
\]
The first-variable condition \(h\in E(p)\) holds in every row.
Thus \eqref{eq:relative-criteria} shows that \(\beta\) fails in both \(U_\alpha\) and
\(U_\eta\), while \(\eta\) fails in both \(U_\alpha\) and \(U_\beta\).
This proves incomparability of \(U_\beta,U_\eta\).
The rule \(\delta\) fails in each of these two varieties,
\(\gamma\) fails in \(U_\delta\), and \(\sigma\) fails in \(U_\gamma\).
All relations in \eqref{eq:covers} are therefore strict.
Since the seven candidates exhaust the interval, they are covering relations.
\end{proof}

\begin{corollary}[finite bases]\label{cor:finitebasis}
Relative to \(\AI\), the seven varieties have the finite bases
\begin{equation}\label{eq:finite-bases}
\B,\qquad \B\cup\{\theta\}\quad
 (\theta\in\{\alpha,\beta,\eta,\delta,\gamma,\sigma\}),
\end{equation}
where \(U_\sigma=A\).  Every displayed basis uses at most three variables.
\end{corollary}

\begin{proof}
Combine Proposition~\ref{prop:upperbasis} with
Theorem~\ref{thm:lattice}.
The standard axioms for \(\AI\) also use at most three variables.
\end{proof}

\begin{proposition}[free-algebra counts]\label{prop:freecounts}
For \(n\geq1\), let \(F_W(n)\) be the free \(n\)-generated algebra
in a variety \(W\) in the seven-element interval.
Let \(d_m\) be the number of antichains, including the empty antichain,
in the subset lattice of an \(m\)-element set.  Then
\begin{equation}\label{eq:free-counts}
|F_U(n)|=(2d_{n-1})^n-1,\qquad |F_A(n)|=4^n-1.
\end{equation}
In particular,
\[
 |F_U(3)|=1727,\qquad |F_A(3)|=63.
\]
\end{proposition}

\begin{proof}
For a fixed first variable \(h\), the contents of long words all contain
\(h\).  Removing \(h\), the inclusion-minimal contents that remain form
an antichain in the subset lattice on the other \(n-1\) variables.
The upper criterion \eqref{eq:upper-criterion} shows that this antichain uniquely determines
the long-word contribution with first variable \(h\), up to the
identities of \(U\); the one-letter summand \(h\) is an independent
choice.  Thus there are \(2d_{n-1}\) choices for each \(h\).
The criteria distinguish all these choices, and every choice is
represented by a polynomial unless no summand at all is selected.
This proves the first count.  In \(A\), only the independent choices
of membership in \(L\) and in \(E\) remain, giving \(4^n-1\).
For \(n=3\), the two-element subset lattice has six antichains,
so \eqref{eq:free-counts} gives \(12^3-1=1727\) and \(4^3-1=63\).
\end{proof}

\begin{corollary}[finite generation by a three-generated algebra]\label{cor:fingen}
For every \(W\in I(S)\),
\begin{equation}\label{eq:three-generated}
W=\V(F_W(3)),\qquad |F_W(3)|\leq1727.
\end{equation}
\end{corollary}

\begin{proof}
The algebra \(F_W(3)\) is a quotient of \(F_U(3)\), hence finite.  The subsets \(\{b\}\) and \(\{c\}\) generate \(\Pplus(S_*)\), since \(\{c\}^2=\{a\}\) and addition is union.  Since \(\Pplus(S_*)\in A\subseteq W\), freeness makes it a homomorphic image of \(F_W(3)\).  Lemma~\ref{lem:endpoint} therefore gives \(A\subseteq\V(F_W(3))\subseteq W\subseteq U\).  By freeness, \(F_W(3)\) and \(W\) have the same identities in at most three variables.  The rules in \eqref{eq:relative-rules} separate all seven interval elements, so no two distinct nodes can occur between \(A\) and \(U\) with the same three-variable theory.  Thus \(\V(F_W(3))=W\).
\end{proof}

\section{The degenerate cases}\label{sec:boundary}

Let \(\mathbf 2=(\{0,1\},\vee,\wedge)\) and \(\Dvar=\V(\mathbf 2)\).  In \(\mathbf 2\), a word is the meet of its variables, and hence
\begin{equation}\label{eq:distributive-criterion}
\Dvar\models u\leq p
 \quad\Longleftrightarrow\quad
 \exists v\in p\; c(v)\subseteq c(u).
\end{equation}
This criterion follows by assigning variables in \(c(u)\) the value one and all other variables zero.

\subsection{The one-element semigroup}

If \(|S|=1\), then \(\Pplus(S)\) is the one-element semiring and \(\mathcal P(S)\cong\mathbf 2\).  Every nontrivial distributive lattice contains a two-element sublattice: from two distinct elements, take their meet if necessary and retain a comparable pair.  Thus every nontrivial subvariety of \(\Dvar\) contains \(\mathbf 2\), and the interval has exactly two elements.

\subsection{Nontrivial left-zero semigroups}

Suppose \(f=\operatorname{id}_S\) and \(|S|\geq2\).  Write
\(A_L=\V(\Pplus(S))\) and \(U_L=\V(\mathcal P(S))\).
Every word has the value of its first variable when its content variables
are nonempty, and is empty otherwise.  Consequently,
\begin{equation}\label{eq:left-zero-criteria}
\begin{aligned}
 A_L\models u\leq p
 &\Longleftrightarrow h(u)\in\{h(v):v\in p\},\\
 U_L\models u\leq p
 &\Longleftrightarrow \exists v\in p\;[h(v)=h(u),\ c(v)\subseteq c(u)] .
\end{aligned}
\end{equation}
For necessity of the second criterion, take distinct \(a,b\in S\);
assign the first variable of \(u\) the set \(\{b\}\), the other
content variables \(\{a\}\), and outside variables the empty set.
If no covering word exists, the right side avoids \(b\).
For the lower criterion use nonempty singleton assignments throughout.
The sufficiency follows from the word evaluation.
In particular \(A_L=\Mod_{\rm ai}\{xy\approx x\}\):
all words reduce to their first variables, and the criterion distinguishes
the remaining sums.
The upper endpoint has the relative basis
\begin{equation}\label{eq:left-zero-basis}
x^2\approx x,\qquad xyz\approx xzy,\qquad xy\leq x.
\end{equation}
Put \(\rho_L:xz\leq xy+z\).  Then
\begin{equation}\label{eq:left-zero-interval}
\begin{gathered}
 I(S)=\{A_L,\ A_L\vee\Dvar,\ U_L\},\\
 A_L\vee\Dvar=U_L\cap\Mod\{\rho_L\}.
\end{gathered}
\end{equation}

To justify the basis \eqref{eq:left-zero-basis}, tail permutation and idempotence reduce
each word to its first variable and content.
The inequality \(xy\leq x\) absorbs larger content into smaller
content with the same first variable.
The upper criterion in \eqref{eq:left-zero-criteria} proves completeness.

We now prove \eqref{eq:left-zero-interval}.
The subsemiring \(\{\varnothing,\{a\}\}\) gives \(\Dvar\subseteq U_L\),
and \(\rho_L\) holds in both \(A_L\) and \(\Dvar\).
If an absorption \(u\leq p\) holds in these two varieties, choose a
right-hand word \(g\) starting with \(h=h(u)\) and a right-hand word
\(v\) with \(c(v)\subseteq c(u)\), using \eqref{eq:left-zero-criteria} and \eqref{eq:distributive-criterion}.
If \(g=h\), the upper criterion already proves the absorption.
Otherwise write \(g=ht\) and use \(\rho_L\) with \(x=h,y=t,z=v\):
\[
 u\leq hv\leq ht+v=g+v\leq p.
\]
It follows that all identities of \(A_L\vee\Dvar\) follow from
\eqref{eq:left-zero-basis} and \(\rho_L\), proving the equality in \eqref{eq:left-zero-interval}.

For exhaustiveness, take any strict absorption valid in \(A_L\) but
not \(U_L\).  Put \(h=h(u)\) and \(C=c(u)\).
Map \(h\) to \(x\), variables in \(C\setminus\{h\}\) to \(z\), and
outside variables to \(xy\).
A right-hand word starting with \(h\) must be external, so its image
is absorbed by \(xy\).
Words starting with other internal variables are absorbed by \(z\);
those starting with outside variables are absorbed by \(xy\).
The left side becomes \(xz\) if \(C\neq\{h\}\), and \(x\) otherwise.
In the latter case \(xz\leq x\) still yields \(\rho_L\).
Thus every proper subvariety of \(U_L\) containing \(A_L\) is contained
in \(A_L\vee\Dvar\).

If an absorption valid in \(A_L\) fails in \(A_L\vee\Dvar\), it fails
in \(\Dvar\); by \eqref{eq:distributive-criterion} every right-hand word is external.
Mapping internal variables to \(x\) and external variables to \(xy\)
gives \(x\leq xy\).
Together with \eqref{eq:left-zero-basis} this yields \(xy\approx x\), and hence \(A_L\).
Any proper subvariety of \(A_L\vee\Dvar\) containing \(A_L\) has such
a separating absorption, so must equal \(A_L\).
Finally \(\Dvar\not\subseteq A_L\), because \(\mathbf2\) fails
\(xy\approx x\); and \(\rho_L\) fails in \(U_L\) at
\(x=\{a\},y=\varnothing,z=\{b\}\).
This proves the strict three-element chain.

\subsection{Nontrivial null semigroups}

Suppose \(|f(S)|=1\), write \(o\) for the common product, and assume \(|S|\geq2\).  Let \(A_N,U_N\) denote the two endpoints.  A long word evaluates to \(\{o\}\) when all its variables are nonempty and to the empty set otherwise.  Hence
\begin{equation}\label{eq:null-criteria}
\begin{aligned}
 A_N\models u\leq p
 &\Longleftrightarrow p\text{ has a long summand},\\
 U_N\models u\leq p
 &\Longleftrightarrow \exists v\in p\;[|v|\geq2,\ c(v)\subseteq c(u)]
\end{aligned}
\end{equation}
for every long word \(u\).
For both endpoints a one-letter absorption \(x\leq p\) holds exactly
when \(x\in L(p)\).
To see its necessity, choose \(a\neq o\), assign \(x\) the value
\(\{a\}\), and all other variables \(\{o\}\).
For necessity in \eqref{eq:null-criteria}, assign every internal variable \(\{a\}\);
use \(\{a\}\) also outside \(c(u)\) for the lower endpoint, and the
empty set outside \(c(u)\) for the upper endpoint.
When the corresponding condition fails, the right side avoids \(o\).
The upper endpoint has the relative basis
\begin{equation}\label{eq:null-basis}
xy\approx yx,\quad x^2\approx x^3,\quad
 xy\approx x^2y^2,\quad x^2y\leq x^2.
\end{equation}
Let \(\rho_N:x^2\leq x+y^2\).  Then
\begin{equation}\label{eq:null-interval}
\begin{gathered}
 I(S)=\{A_N,\ A_N\vee\Dvar,\ U_N\},\\
 A_N\vee\Dvar=U_N\cap\Mod\{\rho_N\}.
\end{gathered}
\end{equation}

For completeness of \eqref{eq:null-basis}, commutativity and the third identity let us
square all letters in a long word; the second identity removes higher
powers.  Thus a long word with content \(C\) reduces to
\(\prod_{x\in C}x^2\).
The last identity absorbs larger content into smaller content.
These reductions derive exactly \eqref{eq:null-criteria}, and the one-letter criterion
then proves completeness.

The subsemiring \(\{\varnothing,\{o\}\}\) gives \(\Dvar\subseteq U_N\),
and \(\rho_N\) holds in both \(A_N\) and \(\Dvar\).
If a long absorption holds in both, \eqref{eq:null-criteria} supplies a long right-hand
word and \eqref{eq:distributive-criterion} supplies an internal right-hand word.
An internal long word already gives the upper criterion.
Otherwise choose an internal one-letter summand \(k\) and a long
right-hand summand \(v\).  Substitution \(x=k,y=v\) in \(\rho_N\) gives
\[
 u\leq k^2\leq k+v^2\approx k+v\leq p.
\]
The one-letter absorptions are already upper identities, so this proves
the relative basis in \eqref{eq:null-interval}.

For a strict absorption valid in \(A_N\) but not \(U_N\), the left side
must be long and there is no internal long right-hand word.
Map variables in its content to \(x\) and all outside variables to
\(y^2\).  The left side becomes \(x^2\).
Internal one-letter summands become \(x\), and every external summand
is absorbed by \(y^2\).
Thus the absorption implies \(\rho_N\), so every proper intermediate
variety is contained in \(A_N\vee\Dvar\).

If an absorption valid in \(A_N\) fails in \(A_N\vee\Dvar\), it
fails in \(\Dvar\), and therefore has no internal right-hand word.
The same substitution gives \(x^2\leq y^2\).
Swapping variables gives equality of all squares.
By \eqref{eq:null-basis}, for any \(x,y\) we have
\(xy\approx x^2y^2\approx x^4\approx x^2\);
all long products consequently coincide.
Every identity of \(A_N\) now follows from \eqref{eq:null-criteria} and the one-letter
criterion, so this additional identity defines \(A_N\) over \(U_N\).
This excludes any further variety below \(A_N\vee\Dvar\).

Finally \(\Dvar\not\subseteq A_N\), since \(x^2\approx y^2\) holds
in \(A_N\) but fails in \(\mathbf2\); and \(\rho_N\) fails in \(U_N\)
at \(x=\{a\},y=\varnothing\).
This proves \eqref{eq:null-interval} and completes the proof of
Theorem~\ref{thm:main}.

\section{Two-generated free algebras and basis ranks}\label{sec:ranks}

Throughout this section assume \eqref{eq:nondegeneracy}.
Let \(L\) and \(N\) be a two-element left-zero semigroup and a
two-element null semigroup, respectively.  Write
\(U_L=\V(\mathcal P(L))\) and \(U_N=\V(\mathcal P(N))\),
and retain \(\Dvar=\V(\mathbf2)\) from Section~\ref{sec:boundary}.

\begin{proposition}[join descriptions]\label{prop:joins}
The intermediate varieties admit the descriptions
\begin{equation}\label{eq:join-descriptions}
\begin{aligned}
 U_\gamma&=A\vee\Dvar, & U_\beta&=A\vee U_N,\\
 U_\eta&=A\vee U_L, & U_\alpha&=A\vee U_N\vee U_L,\\
 U_\delta&=(A\vee U_N)\cap(A\vee U_L).
\end{aligned}
\end{equation}
\end{proposition}

\begin{proof}
An identity holds in a join of varieties if and only if it holds in
each of them.  For a long word \(u=w_{h,C}\), the lower endpoint
requires \(h\in E(p)\) for the absorption \(u\leq p\).
The additional requirements from \(\Dvar,U_N,U_L\) are, respectively,
\(J_C(p),N_C(p),H_{h,C}(p)\), by
\eqref{eq:distributive-criterion}, \eqref{eq:null-criteria}, and
\eqref{eq:left-zero-criteria}.
These are exactly the rows for \(U_\gamma,U_\beta,U_\eta\) in
\eqref{eq:relative-criteria}.
For a one-letter word, the condition from \(A\) already requires
that word to occur as a summand, which suffices in every variety.
This proves the first three descriptions.
The last two follow from Theorem~\ref{thm:lattice}.
\end{proof}

\begin{lemma}[two-variable theories]\label{lem:two-variable}
There are natural isomorphisms
\begin{equation}\label{eq:two-free-isomorphisms}
\begin{aligned}
 F_U(2)&\cong F_{U_\alpha}(2)\cong F_{U_\beta}(2),\\
 F_{U_\eta}(2)&\cong F_{U_\delta}(2)\cong F_{U_\gamma}(2).
\end{aligned}
\end{equation}
The orders are \(35\) and \(24\), respectively, while \(|F_A(2)|=15\).
\end{lemma}

\begin{proof}
Consider an absorption \(u\leq p\) in at most two variables that
holds in \(A\) but fails in \(U\).  By Lemma~\ref{lem:endpoint},
\(u\) is long and some long summand of \(p\) starts with
\(h=h(u)\).
If \(c(u)\) contained every variable occurring in the absorption,
that summand would give the upper criterion, a contradiction.
Thus \(c(u)=\{h\}\).
There is no internal long right-hand word, and the predicates
\(H,T,J\) all reduce to the condition \(h\in L(p)\).
Theorem~\ref{thm:sixtypes} therefore assigns such an absorption to
\(\eta\) or \(\sigma\), according as that condition holds or fails.
The varieties \(U,U_\alpha,U_\beta\) satisfy neither of these rules,
and \(U_\eta,U_\delta,U_\gamma\) satisfy \(\eta\) but not \(\sigma\).
They consequently have the same identities in at most two variables
within each of these two groups.  Their free algebras on two named
generators have identical defining congruences, proving
\eqref{eq:two-free-isomorphisms}.

The values \(35\) and \(15\) follow from
Proposition~\ref{prop:freecounts}, since \(d_1=3\).
To obtain \(24\), write the two generators as \(h,k\).
The six possible contributions with first variable \(h\) in an
upper normal form are: no summand, \(h\), \(h^2\), \(hk\),
\(h+h^2\), and \(h+hk\).
The rule \(\eta\) identifies precisely the last two contributions.
Indeed \(hk\leq h^2\) holds in \(U\), and
\(h^2\leq h+hk\) gives the reverse absorption after adding \(h\).
The criteria in Theorem~\ref{thm:relativecriteria} distinguish the
other choices independently for each first variable: \(L,E\)
distinguish all but \(h^2\) and \(hk\), and these remain distinct
when \(h\) is absent as a one-letter summand.
There are therefore five choices for each of the two first variables.
Excluding the choice with no summand at all gives \(5^2-1=24\).
\end{proof}

For a variety \(W\) generated by one of its finite-rank free algebras,
define its \emph{basis rank} by
\[
 \operatorname{br}(W)=\min\{n\geq1:W=\V(F_W(n))\}.
\]
This measures free generation of a variety; it is distinct from
the number of variables required in an identity basis.

\begin{theorem}[exact basis ranks]\label{thm:ranks}
The two-generated free algebras and basis ranks of the seven interval
elements are as follows:
\begin{equation}\label{eq:rank-table}
\begin{array}{c|r|c|c}
 W&|F_W(2)|&\V(F_W(2))&\operatorname{br}(W)\\ \hline
 U&35&U_\beta&3\\
 U_\alpha&35&U_\beta&3\\
 U_\beta&35&U_\beta&2\\
 U_\eta&24&U_\gamma&3\\
 U_\delta&24&U_\gamma&3\\
 U_\gamma&24&U_\gamma&2\\
 A&15&A&2
\end{array}.
\end{equation}
\end{theorem}

\begin{proof}
As observed in the proof of Corollary~\ref{cor:fingen}, the semiring
\(\Pplus(S_*)\) is generated by \(\{b\},\{c\}\).
Consequently \(A\) is generated by a two-generated algebra.
Write \(N=\{o,d\}\), where every product equals \(o\).
The semiring \(\mathcal P(N)\) is generated by
\(\varnothing,\{d\}\), since \(\{d\}^2=\{o\}\).
The semiring \(\mathbf2\) is generated by its two elements.

Every two-generated member of \(W\) is a homomorphic image of
\(F_W(2)\).  The descriptions of \(U_\beta\) and \(U_\gamma\) in
Proposition~\ref{prop:joins} therefore give
\[
 \V(F_A(2))=A,\qquad
 \V(F_{U_\beta}(2))=U_\beta,\qquad
 \V(F_{U_\gamma}(2))=U_\gamma.
\]
Lemma~\ref{lem:two-variable} now proves the two middle columns of
\eqref{eq:rank-table} for all seven varieties.

For completeness, basis rank one is impossible here.
Proposition~\ref{prop:freecounts} gives
\(|F_U(1)|=|F_A(1)|=3\).
Every \(F_W(1)\) between them is thus the same three-element
algebra, represented by \(x,x^2,x+x^2\); every product in this
algebra equals \(x^2\).
In particular it is commutative, while \(A\subseteq W\) contains
the noncommutative semiring \(\Pplus(S_*)\).
Thus \(\operatorname{br}(W)\geq2\).
The strict inclusions in Theorem~\ref{thm:lattice} and
\eqref{eq:rank-table} exclude rank two for exactly four nodes.
Corollary~\ref{cor:fingen} supplies rank at most three for each,
completing the table.
\end{proof}

\begin{samepage}
\section{The opposite semigroup}

The classification also applies to semigroups satisfying
\(xy\approx y^2\).  If \(S^{\mathrm{op}}\) denotes the opposite
semigroup, then
\[
 \Pplus(S^{\mathrm{op}})\cong\Pplus(S)^{\mathrm{op}},\qquad
 \mathcal P(S^{\mathrm{op}})\cong\mathcal P(S)^{\mathrm{op}}.
\]
Reversing multiplication induces a lattice isomorphism
\(I(S)\cong I(S^{\mathrm{op}})\).  Theorem~\ref{thm:main} therefore
gives the same \(2,3,7\) classification, with right-zero semigroups
in place of left-zero semigroups.
\end{samepage}
\par\penalty0
\begin{samepage}
\section*{Acknowledgements}
\begingroup\small
This work was supported by the Science and Technology Research Program
of the Chongqing Municipal Education Commission
(Grant No.~\mbox{KJZD-K2024011102}) and the Joint Fund for Innovation and
Development of the Chongqing Natural Science Foundation
(Chongqing Municipal Education Commission)
(Grant No.~\mbox{CSTB2025NSCQ-LZX0067}).
\par\endgroup
\end{samepage}

\bibliographystyle{amsplain}
\bibliography{seven-interval}

@book{BurrisSankappanavar,
  author = {Burris, Stanley and Sankappanavar, H. P.},
  title = {A Course in Universal Algebra},
  series = {Graduate Texts in Mathematics},
  volume = {78},
  publisher = {Springer-Verlag},
  address = {New York},
  year = {1981}
}

@misc{GaoRenZhao2026,
  author = {Gao, Zidong and Ren, Miaomiao and Zhao, Xianzhong},
  title = {A continuum in the lattice of semiring varieties: the interval {$[\mathsf{V}(S),\mathsf{V}(S^0)]$}},
  year = {2026},
  howpublished = {Preprint},
  note = {\href{https://arxiv.org/abs/2609.24514}{arXiv:2609.24514}}
}

@article{KurilPolak2005,
  author = {Ku{\v r}il, Martin and Pol{\'a}k, Libor},
  title = {On varieties of semilattice-ordered semigroups},
  journal = {Semigroup Forum},
  volume = {71},
  number = {1},
  pages = {27--48},
  year = {2005}
}

@article{Kuril2017,
  author = {Ku{\v r}il, Martin},
  title = {Admissible closure operators and varieties of semilattice-ordered normal bands},
  journal = {Acta Sci. Math. (Szeged)},
  volume = {83},
  number = {1--2},
  pages = {35--50},
  year = {2017},
  note = {\href{https://doi.org/10.14232/actasm-016-777-4}{doi:10.14232/actasm-016-777-4}}
}

@article{RenZhao2016,
  author = {Ren, Miaomiao and Zhao, Xianzhong},
  title = {The varieties of semilattice-ordered semigroups satisfying {$x^3\approx x$} and {$xy\approx yx$}},
  journal = {Period. Math. Hungar.},
  volume = {72},
  number = {2},
  pages = {158--170},
  year = {2016}
}

@misc{ShaoGao2026,
  author = {Shao, Xiaolei and Gao, Zidong},
  title = {Nonfinitely based intervals of power semiring varieties},
  year = {2026},
  howpublished = {Preprint},
  note = {\href{https://arxiv.org/abs/2609.26019}{arXiv:2609.26019}}
}

@article{WuRenZhao2024,
  author = {Wu, Y. N. and Ren, Miaomiao and Zhao, Xianzhong},
  title = {The additively idempotent semiring {$S_7^0$} is nonfinitely based},
  journal = {Semigroup Forum},
  volume = {108},
  number = {2},
  pages = {479--487},
  year = {2024}
}

@misc{YueRen2025,
  author = {Yue, Mengya and Ren, Miaomiao},
  title = {Every additively idempotent semiring satisfying {$xy\approx xz$} is finitely based},
  year = {2025},
  howpublished = {Preprint},
  note = {\href{https://arxiv.org/abs/2509.16711}{arXiv:2509.16711}}
}

@misc{YueRenGao2026,
  author = {Yue, Mengya and Ren, Miaomiao and Gao, Zidong},
  title = {Explicit equational bases for the power semirings of {$S_7$}},
  year = {2026},
  howpublished = {Preprint},
  note = {\href{https://arxiv.org/abs/2609.19957}{arXiv:2609.19957}}
}

@article{YueRenZengShao2025,
  author = {Yue, Mengya and Ren, Miaomiao and Zeng, L. L. and Shao, Y.},
  title = {The finite basis problem for additively idempotent semirings of order four, {II}},
  journal = {Algebra Universalis},
  year = {2025},
  note = {\href{https://doi.org/10.1007/s00012-025-00908-5}{doi:10.1007/s00012-025-00908-5}}
}

\end{document}